\documentclass[graybox]{svmult}

\usepackage{bbold,amsmath,amssymb,amscd,ntheorem,enumerate,url,tikz}

\newcommand{\bC}{\mathbb{C}}             
\newcommand{\bF}{\mathbb{F}}
\newcommand{\bN}{\mathbb{N}}
\newcommand{\bR}{\mathbb{R}}
\newcommand{\cA}{\mathcal{A}}            
\newcommand{\cB}{\mathcal{B}}
\newcommand{\cF}{\mathcal{F}}
\newcommand{\cG}{\mathcal{G}}
\newcommand{\cH}{\mathcal{H}}
\newcommand{\cL}{\mathcal{L}}
\newcommand{\cP}{\mathcal{P}}
\newcommand{\cS}{\mathcal{S}}
\newcommand{\tM}{\widetilde{M}}          

\newcommand{\tcA}{\widetilde{\cA}}       
\newcommand{\tcS}{\widetilde{\cS}}
\newcommand{\bbone}{\text{\textbb 1}}    
\newcommand{\bone}{\mathbb I}            
\newcommand{\tr}{\operatorname{Tr}}      
\newcommand{\ltwo}{\ensuremath{\ell_{2}}}
\newcommand{\Ltwo}{\ensuremath{L_{2}}}
\newcommand{\norm}[1]{\lVert#1\rVert}

\title{Positive Operator Valued Measures: \\A General Setting for Frames}
\titlerunning{A General Setting for Frames}
\author{Bill Moran, Stephen Howard, and Doug Cochran}
\authorrunning{Moran, Howard, and Cochran} 
\institute{Bill Moran \at Defence Science Institute, University of Melbourne, Parkville VIC, Australia, \email{wmoran@unimelb.edu.au}
\and Stephen Howard \at Defence Science and Technology Organisation, Salisbury SA, Australia \email{stephen.howard@dsto.defence.gov.au}
\and Doug Cochran \at Arizona State University, Tempe AZ, USA, \email{cochran@asu.edu}}

\begin{document}
\maketitle

\abstract{ This paper presents an overview of close parallels that exist
  between the theory of positive operator-valued measures (POVMs) associated
  with a separable Hilbert space and the theory of frames on that space,
  including its most important generalizations.  The concept of a framed POVM
  is introduced, and classical frames, fusion frames, generalized frames, and
  other variants of frames are all shown to to arise as framed POVMs. This
  observation allows drawing on a rich existing theory of POVMs to provide new
  perspectives in the study of frames.  }

\section{Introduction}
\label{sec:intro}

Frames have become a standard tool in signal processing, allowing uniform
description of many linear but non-orthogonal transform techniques that
underpin a wide variety of signal and image processing algorithms.  Initially
popularized in connection with wavelet applications, frames are now a standard
tool in sampling, compression, array processing, as well as in spectral and
other transform methods for time series.

Frames were initially introduced in a 1952 paper of Duffin and Schaeffer
~\cite{Duffin52}, where they appeared as an abstraction of sampled Fourier
transforms. Little interest was shown in them until the appearance of the 1986
paper \cite{Daubechies86} by Daubechies, Grossmann, and Meyer which coincided
with the rise of wavelet methods in signal processing. Subsequently they were
taken up by numerous authors. Several excellent sources, including
\cite{Christensen03,Groechenig01,Feichtinger98,Heil89}, are available for
further details of both the theory and the many applications of frames.

The standard definition of a frame is as a collection $\cF=\{\varphi_k : k\in
K\}$ of elements of a separable Hilbert space $\cH$.  The index set $K$ may be
finite or infinite. In order for $\cF$ to constitute a frame, there must exist
constants $0<A\leq B<\infty$ such that, for all $f\in\cH$,
\begin{equation}
\label{eq:frame_defn}
A\norm{f}^2\leq \sum_{k\in K}\lvert\langle \varphi_k,f\rangle\rvert^2 \leq B\norm{f}^2. 
\end{equation}
Roughly speaking, a projection $f\mapsto \langle f,\varphi_k\rangle$ of a
vector $f$ representing the state of a system onto an individual element
$\varphi_k$ of a frame may be seen as a measurement of that system, and the
aim is to reconstruct the state $f$ from the collection of all individual
measurements $\{ \langle f,\varphi_k\rangle : k\in K\}$ in a robust way. The
frame condition as stated in \eqref{eq:frame_defn} expresses the ability to do
that, and the frame bounds $A$ and $B$ provide a measure of robustness. If
$A=B$, the frame is said to be \emph{tight}. Orthonormal bases are special
cases of tight frames, and for these $A=B=1$.

Several generalizations of the basic concept of a frame have been
proposed. These include, in particular, the possibility that the family
$\{\varphi_k : k\in K\}\subset\cH$ is indexed by a continuum rather than a
discrete index set, resulting in what are called \emph{generalized
  frames}. There are various formulations of generalized frames in the
literature; see in particular \cite{Ali93}. From the perspective of this
paper, the infrastructure of a generalized frame is a measurable function from
a measure space, which serves the role of the index set, to
$\cH$. Specifically, let $(\Omega,\cB,\mu)$ be a measure space (e.g.,
$\Omega=\bR$ with $\cB$ its Borel sets and $\mu$ Lebesgue measure) and let
$\Phi:\Omega\to\cH$ be a $\mu$-measurable function. The collection $\{\Phi(t)
: t\in\Omega\}\subset\cH$ is a generalized frame for $\cH$ if it satisfies a
condition analogous to the frame condition~\eqref{eq:frame_defn}; i.e., for
all $f\in\cH$,
\begin{equation}
\label{eq:gen_framea}
A\norm{f}^2 \leq\int_\Omega \vert\langle\Phi(t),f\rangle\vert^2\,d\mu(t)\leq B\norm{f}^2.
\end{equation}
Define $\Pi_\varphi:\cH\to\cH$ to be orthogonal projection into the
one-dimensional subspace spanned by the unit-norm element $\varphi\in\cH$;
i.e., $\Pi_\varphi(f)=\left<\varphi,f\right>\varphi$. With this notation,
\eqref{eq:gen_framea} becomes
\begin{equation}
\label{eq:gen_frame}
A\bone \leq \int_\Omega \Pi_{\Phi(t)}\,d\mu(t)\leq B\bone ,
\end{equation}
where $\bone$ denotes the identity operator on $\cH$ and the inequalities mean
that the differences are positive definite operators on $\cH$. The integral in
(\ref{eq:gen_frame}) is in the weak sense; i.e., for a suitable measurable
family of operators $\{S(t) : t\in\Omega\}$ on $\cH$, the integral
$\int_\Omega S(t)\,d\mu(t)$ is defined to be the operator $D$ satisfying
\begin{equation*}
\langle f, D\varphi\rangle= \int_\Omega \langle f, S(t)\varphi\rangle\,d\mu(t) 
\end{equation*}
for $f$ and $\varphi$ in $\cH$.

\emph{Fusion frames} generalize the concept of a frame in a different
direction. They have received considerable recent attention in the
signal processing literature; see, for example,
\cite{Fornasier04,Casazza04,Sun06,Asgari05}. In a fusion frame, the
one-dimensional projections $\Pi_{\varphi_k}$ are replaced by
projections $\Pi_k$ onto potentially higher dimensional closed
subspaces $W_k\subset \cH$. Thus a fusion frame $\cF$ is a family
$\{(W_k,w_k) : k\in K\}$ of closed subspaces of $\cH$ and a
corresponding family of weights $w_k\geq 0$ satisfying the frame
condition
\begin{equation}
\label{eq:fusion_frame_cond}
A\norm{f}^2\leq \sum_{k\in K} w_k^2\norm{\Pi_k(f)}^2 \leq B\norm{f}^2
\end{equation}
for all $f\in\cH$. Some authors have promoted fusion frames as a means
of representing the problem of fusion of multiple measurements in, for
example, a sensor network. In this view, each projection corresponds
to a node of the network, and the fusion frame itself, as its name
suggests, provides the mechanism for fusion of these measurements
centrally.

Not surprisingly, the ideas of generalized frames and fusion frames
can be combined into a composite generalization. A \emph{generalized
  fusion frame} $\cF$ for $\cH$ consists of a pair of measurable
functions $(\Phi,w)$. In this setting, $w:\Omega\to\overline{\bR_+}$
and $\Phi:\Omega\to\cP(\cH)$ where $\cP(\cH)$ denotes the space of
orthogonal projections of any rank (including possibly $\infty$) on
$\cH$, endowed with the weak operator topology.  Measurability of
$\Phi$ is in the weak sense that $t\mapsto \langle\varphi,
\Phi(t)\psi\rangle$ is $\mu$-measurable for each $\varphi$ and $\psi$
in $\cH$.  As part of the definition, it is also required that the
function $t\mapsto \Phi(t)f$ is in $\Ltwo(\Omega,\mu)$ for each
$f\in\cH$.  The frame condition in operator form, as in
\eqref{eq:gen_frame}, becomes
\begin{equation*}
A\bone\leq \int_\Omega w(t)^2\Phi(t)\,d\mu(t) \leq B\bone.
\end{equation*} 

As described in later sections of this paper, this definition of a
generalized fusion frame leads to a concept that is, in effect if not
in formalism, remarkably similar to that of a positive operator-valued
measure (POVM) --- a concept that has been prevalent in the quantum
physics literature for many years. This is hardly unexpected from a
signal processing viewpoint, as the concept of POVM was introduced and
developed in quantum mechanics as a means to represent the most
general form of quantum measurement of a system. Further, connections
between POVMs and frames have been noted frequently in the physics
literature (e.g., \cite{Beukema03,Renes04}), although these
relationships seem to be unmentioned in mathematical work on frames.

The remainder of the paper develops a generalization of the POVM
concept as used in quantum mechanics, which encompasses the theory of
frames --- including all of the generalizations discussed above. Once
generalized fusion frames are accepted, setting the discourse in terms
of POVMs enable the importation of much theory from the quantum
mechanics literature and also brings to light some decompositions that
are not readily apparent from the frame formalism.

A key result used in what follows is the classical theorem of
Naimark~\cite{Naimark43} which, long before frames became popular in
signal processing or POVMs were used in quantum mechanics, formalized
analysis and synthesis in this general context. When applied to the
cases above, Naimark's perspective exactly reproduces those notions.

Subsequent sections describe positive operator valued measures,
introduce the theorem of Naimark, and discuss how POVMs relate to
frames and their generalizations.  In this brief description of the
relationship between POVMs and the generalizations of frames, it will
only be possible to touch on the power of the POVM formalism.

\section{Analysis and Synthesis}
\label{sec:analysis_synthesis}

The various concepts of frame, fusion frame, and generalized frame all
give rise to analysis and synthesis operations. In the case of a
frame, a prevalent point of view is that an analysis operator $F$
takes a ``signal'' in $\cH$ to a set of complex ``coefficients'' in
the space $\ltwo(K)$ of square-summable sequences on the index set
$K$; i.e., $F$ is the Bessel map given by $F(f)=\{\langle
f,\varphi_k\rangle : k\in K\}$ where the finiteness of the upper frame
bound $B$ guarantees the square-summability of this coefficient
sequence.  The synthesis operator is the adjoint map
$F^*:\ltwo(K)\to\cH$, given by
\begin{equation*}
F^*(\{a_k\})=\sum_{k\in K} a_k\varphi_k, 
\end{equation*}
and corresponds to synthesis of a signal from a set of
coefficients. It follows directly from \eqref{eq:frame_defn} that the
\emph{frame operator} $\bF=F^*F$ satisfies
\begin{equation}
\label{eq:op_frame_cond}
A\bone\leq \bF\leq B\bone .
\end{equation}

To accommodate developments later in this paper, it is useful to
describe analysis and synthesis with frames in a slightly different
way.  With each $\varphi_k$ in the frame $\cF$, associate the
one-dimensional orthogonal projection operator $\Pi_k$ that takes
$f\in\cH$ to
\begin{equation*}
\Pi_k(f)=\frac{\left<\varphi_k,f\right>}{||\varphi_k||^2}\,\varphi_k
\end{equation*}
Note that $\Pi_k:\cH\to W_k$ where $W_k$ is the one-dimensional subspace of
$\cH$ spanned by $\varphi_k$. Also,
$\norm{\Pi_k(f)}=|\left<\varphi_k,f\right>|/\norm{\varphi_k}$.  Thus the frame
condition (\ref{eq:frame_defn}) is equivalent to
\begin{equation*}
A\norm{f}^2\leq\sum_{k\in K}w_k^2 \norm{\Pi_k(f)}^2 \leq B\norm{f}^2
\end{equation*}
where $w_k=|\left<\varphi_k,f\right>|\geq 0$. From a comparison of this
expression with (\ref{eq:fusion_frame_cond}), it is clear that the weights
$w_k$ account for the possibility that the frame elements $\varphi_k\in\cF$
are not of unit norm. Although it is typical to think of the analysis operator
as producing a set of coefficients for each signal $f\in\cH$ via the Bessel
map, as described above, it is more suitable for generalization to regard it
as a map from $\cH$ to $\cH$ that ``channelizes'' $f$ into signals $w_k
\Pi_k(f)\in W_k\subset\cH$.  The synthesis operator is then a linear rule for
combining a set of signals from the channels $W_k$ to form an aggregate signal
in $\cH$.

With this view, the analysis operator for a fusion frame is a natural
generalization of its frame counterpart in which the subspaces $W_k$ can be of
dimension greater than one and the projection operators $\Pi_k$ are from $\cH$
to $W_k$. The analysis operator is $F:\cH\to \bigoplus_{k\in K} W_k$ given by
\begin{equation*}
F(f)=\{w_k\Pi_k(f) : k\in K\} \in \bigoplus_{k\in K} W_k. 
\end{equation*}
The adjoint map $F^*:\bigoplus_{k\in K} W_k\to\cH$ is given by
\begin{equation*}
F^*(\{\xi_k\})=\sum_{k\in K} w_k \xi_k \in \cH, \quad \{\xi_k : k\in K\}\in \bigoplus_{k\in K} W_k.
\end{equation*}
The frame bound conditions guarantee that everything is well-defined. The
corresponding fusion frame operator $\bF=F^* F:\cH\to\cH$ is given by
\begin{equation*}
   \bF(f)= \sum_{k\in K} w_k^2 \Pi_k(f) ,
\end{equation*}
and the same kind of frame bound inequality as in \eqref{eq:op_frame_cond} holds for fusion frames.

For the generalized frame described in Section \ref{sec:intro}, 
the frame operator $F:\cH\to \Ltwo(\Omega,\mu)$ is given by
\begin{equation*}
F(f)(t)=\langle f, \Phi(t)\rangle, \quad t\in \Omega,\  f\in \cH,
\end{equation*}
and its adjoint by
\begin{equation*}
F^*(u)=\int_\Omega u(t)\Phi(t)\,d\mu(t) \in \cH, \quad u\in
  L^2(\Omega,\mu). 
\end{equation*}
Again, the generalized frame operator $\bF=F^*F$ satisfies inequalities \eqref{eq:op_frame_cond}. 

For generalized fusion frames there is a corresponding definition of analysis
and synthesis operators, but its description requires the definition of direct
integrals of Hilbert spaces \cite{Dixmier81}. In any case the ideas will be
subsumed under the more general development to follow.

It is immediately evident that, in each case discussed above, the synthesis
operator does not reconstruct the analyzed signal; i.e., in general
$F^*F\neq\bone$. In the case of a frame, inversion of the analysis operator is
performed by invoking a \emph{dual frame}. There are various different usages
of this terminology in the literature (see
\cite{Casazza04,Gavruta07,Heil89}). For the purposes here, given a frame
$\{\varphi_k\}$ for the Hilbert space $\cH$, a dual frame
$\{\tilde{\varphi}_k\}$ satisfies
\begin{equation}
  \label{eq:dual_frame}
f=\sum_{k\in K}\langle \varphi_k,f\rangle \tilde{\varphi}_k=\sum_{k\in K}\langle \tilde{\varphi}_k,f\rangle \varphi_k. 
\end{equation}
In other words, the dual frame inverts the analysis and synthesis operations
of the original frame to give perfect reconstruction. Such a dual frame always
exists; indeed, it is easy to verify that
\begin{equation}
  \label{eq:canonical_dual}
\tilde{\varphi}_k=\bF^{-1}(\varphi_k)
\end{equation}
has the appropriate property. Dual frames as defined in \eqref{eq:dual_frame}
are not in general unique; the one in equation~\eqref{eq:canonical_dual} is
called the \emph{canonical dual frame}. In the case of a fusion frame
$\{(W_k,w_k) : k\in K\}$, there also exist dual fusion frames. The canonical
dual fusion frame is $\{(\bF^{-1}W_k,w_k) : k\in K\}$. See \cite{Gavruta07}
for proofs of the existence and discussion of the properties of dual frames in
this context.

\section{Positive Operator Valued Measures}
\label{sec:chap3_posit-oper-valu}

The goal of this section is to define a \emph{framed POVM} and give some
examples of such objects. Consider a topological space $\Omega$ which, to
avoid technicalities, will be assumed to be ``nice;'' e.g., a complete
separable metric space or a locally compact second countable space.  The
crucial point is that $\Omega$ has sufficient structure to make the concept of
regularity of measures meaningful and useful, though regularity will not be
explicitly discussed in this paper. Denote by $\cB(\Omega)$ the
$\sigma$-algebra of Borel sets on $\Omega$ and by $\cP(\cH)$ the space of
positive operators on a Hilbert space $\cH$. A \emph{framed POVM} a function
$M:\cB(\Omega)\to \cP(\cH)$ satisfying the following two conditions:
\begin{enumerate}[POVM-1)]
\item For all $f$ in $\cH$, $\omega\mapsto \langle f,M(\omega)f\rangle$ is a
  regular Borel measure on $\cB(\Omega)$, denoted by $\mu_f$, and
\item $A\bone \leq M(\Omega)\leq B\bone$ for some $0<A\leq B<\infty$. 
\end{enumerate}
As in the case of frames, the numbers $A$ and $B$ are called the frame bounds
for $M$.  Without the condition POVM-2), the object is called a POVM; i.e.,
without the epithet ``framed.''  Such a function is a measure on $\cB(\Omega)$
that takes values in the set of positive operators on $\cH$, though the
countable aspect of its additivity is only in a weak sense.  In the quantum
mechanics context, POVM-2) is replaced by the more strict requirement that
$M(\Omega)=\bone$. A framed POVM is \emph{tight} if $A=B$, and if $A=B=1$, $M$
is a \emph{probability POVM}. Probability POVMs are used in quantum mechanics
as the most general form of quantum measurement.

As an example of a framed POVM, consider a fusion frame $\{(W_k,w_k) : k\in
K\}$ in $\cH$. Define $\Omega=K$ with the $\sigma$-field $\cB(\Omega)$ taken
to be the power set of $\Omega$. Denoting, as above, projection onto $W_k$ by
$\Pi_k$,
\begin{equation}
\label{eq:fus_frame_povm}
M(\omega)=\sum_{k\in \omega}w_k\Pi_k .
\end{equation}
It is straightforward to see that this satisfies both parts of the definition
of a framed POVM, with the frame bounds being the bounds in the definition of
the fusion frame. Thus every fusion frame, and hence every frame, is trivially
represented as a framed POVM.

If $\cF=\{\Phi(t) : t\in\Omega\}$ is a generalized frame for $\cH$, a POVM
$M:\cB(\Omega)\to\cP(\cH)$ can be defined by
\begin{equation}
\label{eq:povm_gen_frame}
M(\omega)=\int_\omega \Pi_{\Phi(t)}\, d\mu(t),
\end{equation}
where $\Pi_{\Phi(t)}$ denotes projection into the one-dimensional subspace of
$\cH$ spanned by $\Phi(t)$. $M$ is a framed POVM with the same frame bounds as
those of $\cF$.

As will be discussed in Section \ref{sec:POVM_analysis_synthesis}, POVMs
provide a rather general framework for analysis and reconstruction of
signals. It will be seen that framed POVMs are only slightly more general than
generalized fusion frames discussed briefly in Section \ref{sec:intro}.  The
impetus for studying POVMs in this context arises in part from the opportunity
to draw on existing theory about POVMs in the physics literature for
development and description of new constructs in signal processing. Some
examples in this paper illustrate this possibility, though much of the
formalism is left for a later paper.

\section{Spectral Measures and the Naimark Theorem}
\label{sec:spectral}

A POVM $S$ is a \emph{spectral} POVM if
\begin{equation*}
  S(\omega_1\cap \omega_2)=S(\omega_1)S(\omega_2),\qquad \omega_1,
\omega_2 \in \cB(\Omega) .
\end{equation*}
Spectral POVMs arise, for example, in the spectral theorem for a Hermitian
operators on Hilbert space (see for example \cite{Sunder97}). If $S$ is a
spectral POVM, then $S(\Omega)$ is a projection, and every $S(\omega)$ with
$\omega\in\cB(\Omega)$ is a projection dominated by $S(\Omega)$; i.e.,
\begin{equation*}
  S(\omega)S(\Omega)=S(\Omega)S(\omega)=S(\omega) .
\end{equation*}
Thus, for any $\omega\in\cB(\Omega$), $S(\omega)$ is completely specified by
its behavior on the closed subspace $S(\Omega)\cH$ of $\cH$. Consequently, for
most purposes it suffices to assume $S(\Omega)=\bone_\cH$. In particular, if a
spectral POVM is framed, then this condition must hold; conversely, imposing
this condition on a spectral POVM ensures that it is framed. Since the
interest here is on framed POVMs, it will be assumed that
$S(\Omega)=\bone_\cH$ whenever a spectral POVM appears in subsequent
discussion in this paper. Note that, while a spectral POVM $S$ need not be
probability POVM in general, the condition that it is framed implies that $S$
will be a probability POVM.  Intuitively, spectral POVMs play an analogous
role relative to framed POVMs to the one played by orthogonal bases relative
to frames; i.e., spectral POVMs generalize orthogonal bases in a sense similar
to that in which framed POVMs generalize frames.

With this machinery in place, it is possible to state the key theorem on POVMs due to
Naimark~\cite{Naimark43}, who formulated the result for POVMs without the
framed condition. The following version is a relatively straightforward adaptation to framed POVMs. 
\begin{theorem}
  \label{thm:stinespring2}
  Suppose $M:\cB(\Omega)\to\cP(\cH)$ is a framed POVM with frame bounds $A$
  and $B$. Then there is an ``auxiliary'' Hilbert space $\cH^\sharp$, a
  spectral POVM $S$ with values in $\cP(\cH^\sharp)$, and a bounded linear map
  $V:\cH^\sharp\to\cH$ such that
\begin{equation*}
M(\omega)=VS(\omega)V^* , \qquad   \omega\in \cB(\Omega)
\end{equation*}
and $A\bone\leq VV^*\leq B\bone$. 
\end{theorem}

For developments later in this chapter, it will be useful to have a sketch of the proof of this theorem.
Given a POVM $M:\cB(\Omega)\to\cP(\cH)$, consider the linear space $\cL$ of $\cH$-valued simple functions on $\Omega$; i.e., finite linear
combinations of functions of the form
\begin{equation}
  \label{eq:simple_h_valued}
  \xi_\omega(t)=
  \begin{cases}
    \xi&\text{ if $t\in \omega$}\\
    0&\text{ otherwise,}
  \end{cases}
\end{equation}
where $\omega\in \cB(\Omega)$ and $\xi\in \cH$. A pre-inner product on $\cL$ is obtained by defining
\begin{equation}
  \label{eq:innerprodL}
  \langle \xi_\omega, \xi'_{\omega'}\rangle_{\cL}=\langle M(\omega)\xi, M(\omega')\xi'\rangle_\cH.
\end{equation}
Completion followed by factoring out zero-length vectors produces $\cH^\sharp$, as a Hilbert space. 
The map from $\cH$ to $\cL$ taking $\xi$ to $ \xi_\Omega$ results in $V^*:\cH\to\cH^\sharp$ and 
$V$ takes $\xi_\omega$ to $M(\Omega)^*M(\omega)\xi$.  The spectral measure $S$ arises first
on $\cL$ as
\begin{equation}
  \label{eq:specmescons}
  S(\omega')(\xi_\omega)=\xi_{\omega\cap\omega'}\qquad \xi\in \cH, \quad
  \omega, \omega'\in \cB(\Omega),
\end{equation}
and then carries over to $\cH^\sharp$.

The collection $(S,\cH^\sharp,V)$ is known as a \emph{Naimark representation}
of the framed POVM $M:\cB(\Omega)\to\cP(\cH)$.  Further, a Naimark
representation is \emph{minimal} if the set
\begin{equation*}
\{S(\omega)V^*\varphi : \varphi\in\cH,\, \omega\in \cB(\Omega)\} 
\end{equation*}
is dense in $\cH^\sharp$. Minimal Naimark representations are essentially
unique in the sense that if $(S,\cH_\sharp,V)$ and $(S',\cH'_\sharp,V')$ are
two such representations for the same $M$, then there is a surjective isometry
$T:\cH_\sharp\to \cH'_\sharp$ such that $V'T=V$ and
$T^{-1}S'(\omega)T=S(\omega)$ for all $\omega\in\cB(\Omega)$. A fashionable
way to handle the Naimark representation in recent literature (see
\cite{Parthasarathy99}) is to convert POVMs to (completely) positive operators
on commutative $C^*$-algebras via integration. In this setting, Naimark's
theorem becomes a special case of Stinespring's theorem \cite{Stinespring55},
which does not require commutativity of the $C^*$-algebra. A full description
of this approach would be tangential to this paper.

\begin{example}
\label{ex:gen_frame_povm_stine}
\emph{Consider a generalized frame $\Phi:\Omega\to\cH$ on the measure space
  $(\Omega,\mu)$ with frame bounds $A\leq B$.  $\Phi$ gives rise to a framed
  POVM $M$ as in \eqref{eq:povm_gen_frame}. To form a Naimark representation
  for $M$, define the Hilbert space $\cH^\sharp$ to be $\Ltwo(\Omega,\mu)$ and
  let the spectral measure $S$ be the canonical one on this space; i.e.,}
\begin{equation*}
S(\omega)f(t)=\bbone_\omega(t)f(t) , \qquad f\in\Ltwo(\Omega,\mu). 
\end{equation*}
\emph{$S$ is clearly a spectral measure since the characteristic functions
  satisfy $\bbone_\omega\bbone_{\omega'}=\bbone_{\omega\cap \omega'}$. The map
  $V:\Ltwo(\Omega,\mu)\to\cH$ is defined by}
\begin{equation*}
V(f)=\int_\Omega f(t)\Phi(t)\, d\mu(t)
\end{equation*}
\emph{where $f(t)\Phi(t)$ is the product of the scalar $f(t)$ and $\Phi(t)\in\cH$. 
It can be verified that this is indeed a (the) minimal Naimark representation of $M$.} 
\end{example}

\begin{example}
\label{ex:fusionframe_povm_stine}
\emph{Let $\cF=\{(W_k, w_k) : k\in K\}$ be a fusion frame in $\cH$.  $\cF$
  corresponds to a framed POVM as in \eqref{eq:fus_frame_povm}.  In this case,
  $\cH^\sharp$ may be taken to be the formal direct sum $\bigoplus_{k\in
    K}W_k$. The appropriate spectral measure $S$ is defined on subsets $J$ of
  $\Omega=K$ by}
\begin{equation}
\label{eq:spect_meas_fusion}
S(J)=\bigoplus_{k\in J} \Pi_k  
\end{equation}
\emph{where $\Pi_k$ is the projection into $W_k$ in $\cH^\sharp$. Writing an
  element $f$ of $\cH^\sharp$ as $f=\{f_k\in W_k : k\in K\}$, the map
  $V:\cH^\sharp\to\cH$ is given by}
\begin{equation}
  \label{eq:fusion_stine_v}
V(f)=\sum_{k\in K} w_k f_k,
\end{equation}
\emph{where the terms in the sum are considered as elements of $\cH$. The
  square-summability of the weights $w_k$ guarantees that the sum in
  \eqref{eq:fusion_stine_v} converges in $\cH$ because the Cauchy-Schwarz
  inequality gives}
\begin{equation}
  \label{eq:sum_frame_fk}
  \sum_{k\in K} \norm{w_k\varphi_k}\leq \Bigl(\sum_{k\in K} w_k^2\Bigr)^{1/2}\Bigl(\sum_{k\in K} \norm{\varphi_k}^2\Bigr)^{1/2} .
\end{equation}
\emph{Thus $V:\cH^\sharp\to \cH$ is a bounded linear map; in fact, by \eqref{eq:sum_frame_fk},}
\begin{equation*}
\norm{V}\leq \Bigl(\sum_{k\in K} w_k^2\Bigr)^{1/2}. 
\end{equation*}
\emph{Its adjoint $V^*:\cH\to\cH^\sharp$ is given by}
\begin{equation*}
V^*(\varphi)=\{w_k\Pi_k(\varphi) : k\in K\}.
\end{equation*}
\emph{Setting $\omega=\Omega=K$ gives $S(\Omega)=\bone$ and
\begin{equation*}
  M(\Omega)=VS(\Omega)V^*=VV^*,
\end{equation*}
The frame bounds imply $A\leq VV^*\leq B$ and, if the fusion frame is tight, then $VV^*=A\bone$.}
\end{example}
From a comparison of the descriptions in Section~\ref{sec:analysis_synthesis}
with the examples given here, it is evident that Naimark's Theorem provides
exactly the machinery for discussing analysis and synthesis operators in a
general context.  This is undertaken in the next section.

\section{Analysis and Synthesis for General POVMs}
\label{sec:POVM_analysis_synthesis}
The preceding examples indicate that the Naimark representation provides a
mechanism for analysis and synthesis in POVMs that precisely extends the
corresponding ideas for frames and fusion frames. To be specific, let
$M:\cB(\Omega)\rightarrow\cP(\cH)$ be a POVM and let $(S,\cH^\sharp,V)$ be the
corresponding mimimal Naimark representation. In this context, $\cH^\sharp$
will be called the \emph{analysis space} and $V^*:\cH\to\cH^\sharp$ the
\emph{analysis operator}. Similarly, $V:\cH^\sharp\to \cH$ will be called the
\emph{synthesis operator}. The use of this terminology is directly analogous
to the way it is used for frames and their generalizations in Section
\ref{sec:analysis_synthesis}. Further, the Naimark representation also
provides a means, via the spectral measure $S$, for keeping track of the
labeling of the POVM.

Analysis of an element $f\in\cH$ is the $\cH^\sharp$-valued measure $\cA$ on $\cB(\Omega)$ defined by
\begin{equation}
\label{eq:povm_analysis}
\cA(f)(\omega)=\hat{f}(\omega)=S(\omega)V^* f \in\cH^\sharp .
\end{equation}
In the case of a frame $\{\varphi_k : k\in K\}$, this measure on subsets of
$\Omega=K$ associates the ``coefficients'' $\left< f,\varphi_k \right> e_k
\in\ltwo(K)$ with the signal $f$, where $\{e_k : k\in K\}$ is the standard
basis of $\ltwo(K)$. Given a measure $\rho:\cB(\Omega)\to\cH^\sharp$ as in
\eqref{eq:povm_analysis}, the synthesis operator takes $\rho$ to
\begin{equation}
\label{eq:povm_synthesis}
\cS(\rho)=V\int_\Omega \,d\rho(t) \in \cH.   
\end{equation}

As the examples in the preceding sections show, these analysis and syntheses
operators correspond precisely to those of classical frames, fusion frames,
and generalized fusion frames.

\section{Isomorphism of POVMs}
\label{sec:isomorphism}

Two POVMs $(M_1, \Omega,\cH_1)$ and $(M_2, \Omega, \cH_2)$ are isomorphic if there is a surjective unitary transformation $U:\cH_1\to \cH_2$ such that $UM_1(\omega)U^{-1}=M_2(\omega)$ for all $\omega\in \cB(\Omega)$.  The following result is a straightforward consequence of the proof of the Naimark theorem.
\begin{theorem}
\label{thm:naimark_functor}
Suppose that POVMs $M_1:\cB(\Omega)\to\cP(\cH_1)$ and $M_2:\cB(\Omega)\to\cP(\cH_2)$ are isomorphic via the unitary transformation $U:\cH_{1}\to \cH_{2}$. Let $(S_1,\cH^\sharp_1,V_1)$ and $(S_2,\cH^\sharp_2,V_2)$ be minimal Naimark representations of $M_1$ and $M_2$,  respectively. Then there is a unitary transformation $U^\sharp: \cH^\sharp_1\to\cH^\sharp_2$ such that $U^\sharp S_1(\omega)(U^\sharp)^{-1}=S_2(\omega)$ for all $\omega\in\cB(\Omega)$ and the following diagram commutes:
\begin{center}
  \begin{tikzpicture}[scale=2.50]%
  \node (H1) {$\cH_{1}$};
  \node [right of=H1] (H2) {$\cH_{2}$};
  \node [ below of=H1] (H1sharp) {$\cH_{1}^\sharp$};
  \node [ below of=H2] (H2sharp) {$\cH_{2}^\sharp$};
  \draw[->] (H1) to node[above] {$U$} (H2);
  \draw[->] (H1sharp) to node[below] {$U^\sharp$} (H2sharp);
  \draw[<-] (H1) to node[left] {$V_{1}$} (H1sharp);
  \draw[<-] (H2) to node[right] {$V_{2}$} (H2sharp); 
\end{tikzpicture}.
\end{center}
\end{theorem}

Although this result does not appear to be explicitly stated in the literature, it is implicit in many applications of the Naimark and Stinespring theorems. In particular, the paper of Arveson~\cite{Arveson69} discusses related ideas. The proof follows by consideration of the construction of the Naimark representation using Hilbert space valued functions as described in Section \ref{sec:spectral}. Specifically, using the notation of the sketch proof of Naimark's theorem given in Section~\ref{sec:spectral}, observe that for the isomorphic POVMS $M_1$ and $M_2$, $U$ gives rise to a map $\cL_1\to\cL_2$ taking $\xi_\omega$ to $U(\xi)\omega$ which then produces $U^\sharp$. Moreover it follows from the definition of the spectral measure in \eqref{eq:specmescons} that $U^\sharp S_1(\omega)(U^\sharp)^{-1}=S_2(\omega)$ for all $\omega\in\cB(\Omega)$.

\section{Canonical Representations and POVMs}
\label{sec:canonical}

Combining the Naimark theorem and Theorem \ref{thm:naimark_functor} with the canonical
representation of spectral POVMs (described in, e.g., \cite{Sunder97})
yields a canonical representation of POVMs such that two
isomorphic POVMs have the same canonical representation. This
serves to illustrate the utility of the POVM formalism.  The canonical
representation decomposes $\cH^\sharp$, the analysis space of a POVM
$M:\cB(\Omega)\to\cP(\cH)$ that arises in its Naimark representation,
into a direct sum $\bigoplus_{n\in\bN}\cG_n$ such that:
\begin{enumerate}
\item Each of the spaces $\cG_n$ is invariant under the spectral measure; i.e.,
\begin{equation*}
S(\omega)\cG_n\subset \cG_n\qquad \omega\in\cB(\Omega), n\in\bN ,
\end{equation*}
and
\item The restriction of $S$ to $\cG_n$ has uniform multiplicity; i.e.,
  $\cG_n\simeq \bC^{u_n}\otimes \Ltwo(\mu_n)$ if $\cG_n$ has finite dimension
  $u_n$, and $\cG_n\simeq \ltwo(\bN)\otimes\Ltwo(\mu_n)$ if $\cG_n$ is
  infinite-dimensional.
\end{enumerate}
This representation is essentially unique up to unitary equivalence and
replacement of each of the measures $\mu_n$ by one having the same null
sets. Denote by $P_n$ the projection into $\cG_n$, regarded as a subspace of
$\cH^\sharp$. Under the (minimal) Naimark representation, $V:\cH^\sharp\to
\cH$ is such that $V^*S(\omega)V=M(\omega)$ for $\omega\in\cB(\Omega)$.  $V$
can be decomposed as $V=\sum_n VP_n=\sum_n V_n$.  The image of $V^*_n$ is in
$\cG_n$, so that
\begin{equation*}
M(\omega)=\sum_n V_n S(\omega)V_n^* , \qquad \omega\in\cB(\Omega).
\end{equation*}
The map $M_n:\cB(\Omega)\to\cP(\cG_n)$, defined by $M_n(\omega)=V_n
S(\omega)V_n^*$, is a POVM; more precisely, the values of $M_n$ are positive
operators on the closure of the image of $V_n$. The individual measures $M_n$
are themselves POVMs, though they need not be framed even if $M$ is
framed. However, $M_n(\Omega)=V_nV_n^*$.  Thus if $M$ is a framed POVM with
frame bounds $A\leq B$,
\begin{equation*}
  A\bone_\cH\leq \sum M_n(\Omega)=\sum_n V_nV_n^*\leq B\bone_\cH. 
\end{equation*}
Observe that $V_n^*V_nV_m^*V_m=0$ for $n\neq m$, since the image of $V_m^*$
lies in $\cG_m$ which is in the kernel of $V_n$. So, an obvious sense,
\begin{equation*}
  M=\sum_{n\in\bN} M_n.
\end{equation*}
Thus every framed POVM is a sum of ``uniform multiplicity'' POVMs,
though these need not be framed, and this composition is essentially
unique. The canonical representation is characterized by the sequence
of equivalence classes of measures $\{[\mu_n] | n\in\bN\}$ and the linear map
$V$.

\begin{example}
\label{ex:frame_multip}
\emph{Consider a frame $\cF=\{\varphi_k : k\in K\}$ in $\cH$ with frame bounds
  $A\leq B$ and its corresponding framed POVM $M$. In this case $\cH^\sharp$
  is $\ltwo(K)$, $V:\cH^\sharp\to\cH$ is given by $V(e_k)=\varphi_k$. The
  spectral measure on the subsets of $K$ is given by}
\begin{equation*}
S(J)=\sum_{k\in J} \Pi_k , \qquad J\subset K,
\end{equation*}
\emph{where $\Pi_k$ denotes projection into the subspace of $\ltwo(K)$ spanned
  by the standard basis element $e_k$. Alternatively this can be redefined by
  regarding members of $\ltwo(K)$ as complex-valued functions on $\Omega=K$
  and taking $S(J)(f)=\bbone_J f$ so that the spectral measure is uniform with
  multiplicity one.}
\end{example}

\begin{example}
\label{ex:fusionframe_povm_stineddd}
\emph{The the case of a fusion frame $\{(W_k,w_k) : k\in K\}$ is more
  complicated than the frame case.  The spectral measure $S$ on subsets of
  $\Omega = K$ is given by \eqref{eq:spect_meas_fusion}. For each $j\in K$,
  denote}
\begin{equation*}
U_j=\{k\in K: \dim W_k = j \} .
\end{equation*}
\emph{Then}
\begin{equation*}
Y_j= \bigoplus_{k\in U_k} W_k \subset \cH^\sharp
\end{equation*}
\emph{is isomorphic to $\bC^j \otimes \ltwo(U_k)$ or, if $j=\infty$,
  $\ltwo(U_j)\otimes \ltwo(U_j)$. Evidently, $Y_j$ has uniform multiplicity
  $j$, and the measure $\mu_j$ is counting measure on $U_j$, provided $U_j$ is
  not empty. If all $W_k$ for $k\in K$ have the same dimension, then the
  spectral measure $S$ has uniform multiplicity.}
\end{example}

\section{Dual POVMs}
\label{sec:dual_povm}

As observed in Section~\ref{sec:analysis_synthesis}, each of frame
generalizations associates a ``dual'' object with the frame, and there is a
canonical dual in each case. This is also possible for framed POVMs, and
indeed is relatively straightforward using the Naimark
representation. Consider a POVM $M:\cB(\Omega)\to\cP(\cH)$ and its minimal
Naimark representation $(\Omega, S, \cH^\sharp,V)$. The \emph{canonical dual
  POVM} to $M$ is the POVM $\tM:\cB(\Omega)\to\cP(\cH)$ having Naimark
representation $(\Omega,S,\cH^\sharp,(VV^*)^{-1}V)$; i.e.,
\begin{equation*}
\tM(\omega)=(VV^*)^{-1}VS(\omega)V^*(VV^*)^{-1}.
\end{equation*}
The frame condition on $M$ guarantees $0<A\leq V^* V \leq B<\infty$, which not
only ensures the existence of $(VV^*)^{-1}$, but implies $\tM$ is a framed
POVM with bounds $B^{-1}\leq A^{-1}$ (see Theorem~\ref{thm:stinespring2}).
Further,
\begin{equation*}
\begin{aligned}
  M(\omega)\tM(\omega)&=\bigl(VS(\omega)V^*\bigr)\bigl(
  (VV^*)^{-1}VS(\omega)V^*(VV^*)^{-1}\bigr)\\ 
  \tM(\omega)M(\omega)&=\bigl((VV^*)^{-1}VS(\omega)V^*(VV^*)^{-1}\bigr)\bigl(
  VS(\omega)V^*\bigr)
 \end{aligned}
\end{equation*}
In particular, invoking the assumption $S(\Omega)=\bone_{\cH^\sharp}$ gives
\begin{equation*}
  M(\Omega)\tM(\Omega)=\tM(\Omega)M(\Omega)=\bone_{\cH}. 
\end{equation*}

From the point of view of analysis and synthesis, if $f\in\cH$, its analysis
with respect to $M$ is the measure $\cA(f)$ given in
\eqref{eq:povm_analysis}. Subsequently applying the synthesis operator $\tcS$
associated with the canonical dual POVM $\tM$ yields \eqref{eq:povm_synthesis}
gives
\begin{equation*}
  \tcS\cA(f)(\Omega)=(VV^*)^{-1}V S(\Omega)V^*f=f.  
\end{equation*}
Similarly, analysis of $f$ by $\tM$ followed by synthesis with $M$ is also the identity; i.e.,
\begin{equation*}
\cS\tcA(f)(\Omega)=VS(\Omega)V^*(VV^*)^{-1}f=f .
\end{equation*}

\section{Radon-Nikodym Theorem for POVMs}
\label{sec:radon-nikodym}

This section summarizes some results pertinent to framed POVMs on
finite-dimensional Hilbert spaces.  This setting is prevalent in signal
processing applications, and it will be seen that the theory developed is
valid in a number of infinite-dimensional examples as well.  In this setting,
the concept of a framed POVM identical to that of a generalized fusion frame,
described in Section \ref{sec:intro}.

Let $M:\cB(\Omega)\to\cP(\cH)$ be a framed POVM where $\dim \cH$ is
finite. The finite-dimensional assumption on $\cH$ allows definition of a
real-valued Borel measure $\mu(\omega)=\tr(M(\omega))$ on the Borel sets of
$\Omega$. This positive regular Borel measure is a key element in the
following \emph{Radon-Nikodym theorem} for POVMs (see \cite{Chiribella10}).
\begin{theorem}
\label{thm:rdthm}
Let $M:\cB(\Omega)\to\cP(\cH)$ be a POVM with $\cH$ finite-dimensional. Then
there exists a regular positive real-valued measure $\mu$ on $\cB(\Omega)$ and
an operator-valued bounded measurable function $r:\Omega\to \cP(\cH)$ such
that
\begin{equation*}
M(\omega)=\int_{\omega}r(t)\,d\mu(t) , \qquad \omega\in \cB(\Omega). 
\end{equation*}
\end{theorem}

The measure $\mu$ is called the \emph{base measure} of the POVM and $r$ the
\emph{Radon-Nikodym derivative} of the POVM $M$ with respect to $\mu$.  This
representation is useful in facilitating constructions of POVMs when $\cH$ is
finite-dimensional.
\begin{corollary}
\label{thm:rn_frame_bnds}
If $M$ is a framed POVM with frame bounds $A\leq B$, then
\begin{equation*}
A\bone_{\cH}\leq \int_\Omega r(t) \, d\mu(t) \leq B\bone_{\cH}
\end{equation*}
\end{corollary}

It is instructive to observe how this Radon-Nikodym theorem manifests in the
motivating examples. In particular, this result shows that, when $\cH$ is
finite-dimensional, framed POVMs correspond exactly to generalized fusion
frames as introduced in Section \ref{sec:intro}.

\begin{example}
\label{ex:frame_rn}
\emph{Let $\cF=\{\varphi_k : k\in K\}$ be a frame in $\cH$. The associated
  POVM is given by $M(J)=\sum_{k\in J} \Pi_k$ for subsets $J$ of $\Omega=K$.
  In this case, the operator-valued function $r$ is given by}
\begin{equation*}
r(k)=\Pi_k , \qquad k\in K. 
\end{equation*}
\emph{In this special case, there is no need for the finite-dimensional
  restriction on $\cH$. A POVM constructed from a frame in this way
  automatically possesses a Radon-Nikodym derivative with respect to counting
  measure on the subsets of $K$.}
\end{example}

\begin{example}
\label{ex:gen_frame_rn}
\emph{In the case of a generalized frame $\Phi:\Omega\to \cH$ for a Hilbert
  space $\cH$, the associated POVM is given in \eqref{eq:povm_gen_frame}. In
  this case, the operator-valued function is $r(t)=\Pi_{\Phi(t)}$.  As in the
  previous case, a POVM constructed in this way satisfies a Radon-Nikodym
  theorem with respect to the given measure $\mu$ on $\Omega$ even when $\cH$
  is not finite-dimensional.}
\end{example}

\begin{example}
\label{ex:fus_frame_rn}
\emph{For a fusion frame $\{(W_k,w_k) : k\in K\}$, $\Omega= K$ and $\mu$ is
  counting measure on subsets of $K$. The function $r:\cB(K)\to \cP(\cH)$ is
  given by $r(k)=w_k^2\Pi_{W_k}$, which coincides with the previous
  observation that the POVM in this case is defined by
\begin{equation*}
M(\omega)=\sum_{k\in \omega} w^2_k\Pi_{W_k} , \qquad \omega\subset K. 
\end{equation*}
Although the values of $r$ are not projections, they are non-negative
multiples of projections.  If the counting measure $\mu$ were replaced by
$\nu(k)=w_k^2$, then the expression \eqref{eq:povm_gen_frame} for $M$ would
become
\begin{equation*}
M(\omega)=\int_\omega \Pi_{W_k}\,d\nu(k) , \qquad \omega\subset K,
\end{equation*}
and the Radon-Nikodym derivative of $M$ with respect to $\nu$ would have true
projections as its values.}
\end{example}

A POVM $M:\Omega\to\cP(\cH)$ is \emph{decomposable} if there is an essentially bounded 
measurable function $r:\Omega\to \cP(\cH)$ and a measure $\mu$ on $\cB(\Omega)$ such that
\begin{equation*}
M(\omega)=\int_\omega r(t)\,d\mu(t) \qquad \omega\in \cB(\Omega). 
\end{equation*}
As observed above, if $\dim\cH$ is finite, the POVM is decomposable. Further,
every POVM arising from a (generalized) frame is decomposable, even when $\cH$
is not finite-dimensional. In effect, decomposable framed POVMs correspond to
generalized fusion frames as described in Section \ref{sec:intro}, and thus
this concept captures the simultaneous generalization of frames to fusion
frames and generalized frames.

\section{Conclusions}

In this overview, we have set forth the concept of a framed positive
operator-valued measure and shown that classical frames, as well as several
generalizations of frames, arise as special cases of framed POVMs. We have
described how Naimark's theorem for POVMs leads to notions of analysis and
synthesis for POVMs that subsume their frame counterparts.  We have further
discussed how canonical representations of spectral POVMs lead to canonical
descriptions of framed POVMs, and that this leads to a notion of a canonical
dual POVM analogous to that of the canonical dual of a frame.

\end{document}